\documentclass{article}
\usepackage{amsfonts}
\usepackage{amsmath}

\newtheorem{theorem}{Theorem}

\newtheorem{proposition}[theorem]{Proposition}
\newtheorem{remark}[theorem]{Remark}

\begin{document}
	\title{Lagrange and Hamilton geometries applied to a dynamical sistem governing COVID-19 disease}
	\author{Ana-Maria Boldeanu and Mircea Neagu}
	\date{}
	\maketitle
	
	\begin{abstract}
		In this paper we develop, via the least squares variational method,
		the Lagrange-Hamilton geometry (in the sense of nonlinear connections, d-torsions
		and Lagrangian Yang-Mills electromagnetic-like energy) produced by a dynamical
		system governing the spreading of COVID-19 disease. The Jacobi stability of this
		dynamical system is also discussed.
	\end{abstract}
	
	\textbf{Mathematics Subject Classification (2020):} 00A69, 37J05, 70H40.
	
	\textbf{Key words and phrases}: COVID-19 dynamical system, (co)tangent bundles, least squares Lagrangian and Hamiltonian, Lagrange-Hamilton geometry.
	
	\section{Introduction}
	
The spreading of COVID-19 disease is governed by the following dynamical system
investigated in the paper \cite{Gumel}:
	
	\begin{equation} \label{1}
		\left\{\begin{array}{l}\medskip
		\dfrac{dS}{dt} = -\dfrac{\beta_{s}I_{s}+\beta_{a}I_{a}+\beta_{h}I_{h}}{N} \cdot S\\\medskip
		\dfrac{dE}{dt} =\dfrac{\beta_{s}I_{s}+\beta_{a}I_{a}+\beta_{h}I_{h}}{N} \cdot S - \sigma E \\\medskip 
		\dfrac{dI_{s}}{dt} = (1-r)\sigma E - (\phi_{s}+\gamma_{s}+\delta_{s})I_{s} \\\medskip
		\dfrac{dI_{a}}{dt} = r\sigma E - \gamma_{a}I_{a} \\\medskip
		\dfrac{dI_{h}}{dt} = \phi_{s}I_{s}-(\gamma_{h}+\delta_{h})I_{h} \\\medskip
		\dfrac{dR}{dt} =  \gamma_{s}I_{s}+\gamma_{a}I_{a}+\gamma_{h}I_{h},
		\end{array}\right.
	\end{equation}%
	where
	\begin{itemize}
		\item [$\rightarrow$] $N(t) = S(t)+E(t)+I_{s}(t)+I_{a}(t)+I_{h}(t)+R(t)$ is the population size at time $\textit{t}$ which is sub-divided into the compartments of \textit{susceptible} (i.e., people who are at risk of acquiring infection, but have not yet contracted the disease) ($S(t)$), \textit{exposed} (i.e., newly-infected individuals who are incubating the disease) ($E(t)$), $\textit{symptomatically-infectious}$ (i.e., infectious people showing
		clinical symptoms of the disease) ($I_{s}(t)$), $\textit{asymptomatically-infectious}$ (i.e.,
		infectious people showing no clinical symptoms of the disease) ($I_{a}(t)$), $\textit{hospitalized}$ ($I_{h}(t)$) and $\textit{recovered}$ ($R(t)$) individuals.
		\item [$\rightarrow$] The parameters $\beta_{s}, \beta_{a}$ and $\beta_{h}$ represent the rate at which symptomatically-infectious,
		asymptomatically-infectious, and hospitalized individuals transmit
		COVID-19 to susceptible individuals.
		\item [$\rightarrow$] Exposed individuals progress out of the \textit{E} class at a rate $\sigma$ (i.e., $1/ \sigma$ is the intrinsic incubation period of COVID-19).
		\item [$\rightarrow$] It is assumed that the parameter $0 < r \leq 1$ of exposed individuals show no clinical symptoms of COVID-19 (and move to the $I_{a}$ compartment) at the end of the incubation period.
		\item [$\rightarrow$] Individuals in the $I_{s}(I_{a})(I_{h})$ compartment recover from COVID-19 infection at a rate $\gamma_{s}(\gamma_{a})(\gamma_{h})$.
		\item [$\rightarrow$] Infectious individuals are hospitalized (or isolated either at home or in hospital) at a rate $\phi_{s}$.
		\item [$\rightarrow$] Individuals in the symptomatically-infectious $(I_{s})$ and hospitalized $(I_{h})$ compartments die of COVID-19 at a rate $\delta_{s}$ and $\delta_{h}$. One consider that the rate of
		change of the population of COVID-19 -- deceased individuals (which is denoted
		by $D(t)$) is given by the formula
		$$\dfrac{dD}{dt}=\delta_{s}I_{s}+\delta_{h}I_{h}.$$
	\end{itemize}
	\par For more details about this mathematical model for COVID-19 disease finds in
	the paper \cite{Gumel}.

	\section{From COVID-19 dynamical system to Lagrange geometry and Jacobi stability problem}
		
	In the context of COVID-19 dynamical system ($\ref{1}$), we can consider that we work
	on the particular 6-dimensional manifold $M = \mathbb{R}^6$, whose coordinates are
	$$(x^1,x^2,x^3,x^4,x^5,x^6)=(S,E,I_{s},I_{a},I_{h},R).$$
	Moreover, its corresponding tangent $TM$ and cotangent $T^{*}M$ bundles have the
	coordinates $(x^i,y^i)_{i=\overline{1,6}}$, respectively $(x^i,p_i)_{i=\overline{1,6}}.$ 
	\par Now, let us take the vector field $X=(X^i(S,E,I_{s},I_{a},I_{h},R))_{i=\overline{1,6}}$ on the manifold $M = \mathbb{R}^6$, which is described by
	\begin{equation} \label{2}
		\begin{array}{l}\medskip
			X^1(S,E,I_{s},I_{a},I_{h},R)= -\dfrac{\beta_{s}I_{s}+\beta_{a}I_{a}+\beta_{h}I_{h}}{N} \cdot S,  \\\medskip
			X^2(S,E,I_{s},I_{a},I_{h},R)=\dfrac{\beta_{s}I_{s}+\beta_{a}I_{a}+\beta_{h}I_{h}}{N} \cdot S - \sigma E,  \\\medskip
			X^3(S,E,I_{s},I_{a},I_{h},R)=(1-r)\sigma E - (\phi_{s}+\gamma_{s}+\delta_{s})I_{s}, \\\medskip 
			X^4(S,E,I_{s},I_{a},I_{h},R)=r\sigma E - \gamma_{a}I_{a}, \\\medskip 
			X^5(S,E,I_{s},I_{a},I_{h},R)=\phi_{s}I_{s}-(\gamma_{h}+\delta_{h})I_{h},  \\\medskip
			X^6(S,E,I_{s},I_{a},I_{h},R)=\gamma_{s}I_{s}+\gamma_{a}I_{a}+\gamma_{h}I_{h}.
		\end{array}
	\end{equation}
	Then, we can regard the COVID-19 dynamical system ($\ref{1}$) as the system
	\begin{equation} \label{3}
		\dfrac{dx^i}{dt} = X^i(x(t)), \; i=\overline{1,6}.
	\end{equation}
	Obviously, the solutions of class $C^2$ of the dynamical system ($\ref{3}$) are the global
	minimum points for the \textit{least squares Lagrangian} $L : TM \rightarrow \mathbb{R}$, which is given by\footnote{%
		The Latin indices $i,j,k,...$ run from $1$ to $6.$ Moreover, the Einstein
		convention of summation is adopted all over this work.}%
	\begin{equation} \label{4}
		\begin{array}{c}
			L(x,y)=\delta _{ij}\left( y^{i}-X^{i}(x)\right) \left( y^{j}-X^{j}(x)\right)
			\Leftrightarrow \medskip \\ 
			L(x,y)=\left( y^{1}-X^{1}(x)\right) ^{2}+\left( y^{2}-X^{2}(x)\right)
			^{2}+\left( y^{3}-X^{3}(x)\right) ^{2}+ \medskip \\
			+\left( y^{4}-X^{4}(x)\right) ^{2}+\left( y^{5}-X^{5}(x)\right) ^{2}+\left( y^{6}-X^{6}(x)\right) ^{2}.%
		\end{array}
	\end{equation}%
	
	The Euler-Lagrange equations of the least squares Lagrangian (\ref{4}) are expressed by%
	\[
	\frac{\partial L}{\partial x^{k}}-\frac{d}{dt}\left( \frac{\partial L}{%
		\partial y^{k}}\right) =0,\quad y^{k}=\frac{dx^{k}}{dt},\quad k=\overline{1,6}. 
	\]%
	These can be rewritten in the following geometrical form:
	\begin{equation} \label{5}
		\frac{d^{2}x^{k}}{dt^{2}}+2G^{k}(x,y)=0,\quad k=\overline{1,6},  
	\end{equation}%
	where the object $G=(G^k)_{k=\overline{1,6}}$ described by
	\[
	G^{k}=\frac{1}{2}\left( \frac{\partial ^{2}L}{\partial x^{j}\partial y^{k}}%
	y^{j}-\frac{\partial L}{\partial x^{k}}\right) =-\frac{1}{2}\left[ \left( 
	\frac{\partial X^{k}}{\partial x^{j}}-\frac{\partial X^{j}}{\partial x^{k}}%
	\right) y^{j}+\frac{\partial X^{j}}{\partial x^{k}}X^{j}\right] 
	\]%
	is endowed with the geometrical meaning of \textit{semispray} of $L$.
	
	It is already an well-known fact that following the geometrical ideas from the works Miron and
	Anastasiei \cite{Mir-An}, Udri\c{s}te and Neagu \cite{Udr}, \cite{Nea-Udr},
	and Balan and Neagu \cite{Bal-Nea}, we can construct an entire natural
	collection of nonzero Lagrangian geometrical objects (such as nonlinear
	connection, distinguished (d-)torsions and Yang-Mills electromagnetic-like energy) that
	characterize the dynamical system (\ref{3}) and, obviously,
	the COVID-19 system (\ref{1}).
	
	The Jacobian matrix $J=\left( \frac{\partial
		X^{i}}{\partial x^{j}}\right) _{i,j=\overline{1,6}}$ of the vector field $X$ is given by:
			$$J=\begin{pmatrix}
				J_{11} & J_{12} & J_{13} &J_{14} &J_{15} & J_{16} \\
				J_{21} & J_{22} & J_{23} &J_{24} &J_{25} & J_{26} \\ 
				J_{31}=0 & J_{32} & J_{33} & J_{34}=0 & J_{35}=0 & J_{36}=0\\ 
				J_{41}=0 & J_{42} & J_{43}=0 & J_{44} & J_{45}=0 & J_{46}=0 \\ 
				J_{51}=0 & J_{52}=0 & J_{53} & J_{54}=0 & J_{55} & J_{56}=0\\ 
				J_{61}=0 & J_{62}=0 & J_{63} & J_{64} & J_{65} & J_{66}=0\\ 
			\end{pmatrix},$$ \\
	where 
	\[
	J_{11}=-\dfrac{\beta_{s}I_{s}+\beta_{a}I_{a}+\beta_{h}I_{h}}{N^2} \cdot S -\dfrac{\beta_{s}I_{s}+\beta_{a}I_{a}+\beta_{h}I_{h}}{N}, 
	\]
	\[
	J_{12}=-\dfrac{\beta_{s}I_{s}+\beta_{a}I_{a}+\beta_{h}I_{h}}{N^2} \cdot S, \;
	J_{13}=\dfrac{-\beta_{s}N+\beta_{s}I_{s}+\beta_{a}I_{a}+\beta_{h}I_{h}}{N^2} \cdot S,
	\]
	\[
	J_{14}=\dfrac{-\beta_{a}N+\beta_{s}I_{s}+\beta_{a}I_{a}+\beta_{h}I_{h}}{N^2} \cdot S, 
	\]
	\[
	J_{15}=\dfrac{-\beta_{h}N+\beta_{s}I_{s}+\beta_{a}I_{a}+\beta_{h}I_{h}}{N^2} \cdot S,
	\]
	\[
	J_{16}=-\dfrac{\beta_{s}I_{s}+\beta_{a}I_{a}+\beta_{h}I_{h}}{N^2} \cdot S,
	\]
	\[
	J_{21}=\dfrac{\beta_{s}I_{s}+\beta_{a}I_{a}+\beta_{h}I_{h}}{N^2} \cdot S +\dfrac{\beta_{s}I_{s}+\beta_{a}I_{a}+\beta_{h}I_{h}}{N}, 
	\]
	\[
	J_{22}=\dfrac{\beta_{s}I_{s}+\beta_{a}I_{a}+\beta_{h}I_{h}}{N^2} \cdot S-\sigma, \;
	J_{23}=\dfrac{\beta_{s}N-(\beta_{s}I_{s}+\beta_{a}I_{a}+\beta_{h}I_{h})}{N^2} \cdot S,
	\]%
	\[
	J_{24}=\dfrac{\beta_{a}N-(\beta_{s}I_{s}+\beta_{a}I_{a}+\beta_{h}I_{h})}{N^2} \cdot S, 
	\]
	\[
	J_{25}=\dfrac{\beta_{h}N-(\beta_{s}I_{s}+\beta_{a}I_{a}+\beta_{h}I_{h})}{N^2} \cdot S,
	\]
	\[
	J_{26}=\dfrac{\beta_{s}I_{s}+\beta_{a}I_{a}+\beta_{h}I_{h}}{N^2} \cdot S,
	\]
	\[
	J_{32}=(1-r)\sigma, \;
	J_{33}=- (\phi_{s}+\gamma_{s}+\delta_{s}), \;
	J_{42}=r\sigma, \;
	J_{44}=-\gamma_{a},
	\]
	\[
	J_{53}=\phi_{s}, \; 
	J_{55}=-(\gamma_{h}+\delta_{h}), \;
	J_{63}=\gamma_{s}, \;
	J_{64}=\gamma_{a}, \;
	J_{65}=\gamma_{h}.
	\]
	
	As a consequence, the entries of the canonical nonlinear connection $
	N=\left( N_{j}^{i}\right) _{i,j=\overline{1,6}} 
	$ are defined by the formulas $%
	N_{j}^{i}=\partial G^{i}/\partial y^{j}$ (see \cite{Mir-An}), which imply the
	following matrix relation (see \cite{Bal-Nea})%
	\[
	N=-\displaystyle{\frac{1}{2}}\left[ J-J^{t}\right] . 
	\]%
	
	\begin{proposition}
		The \textbf{Lagrangian canonical nonlinear connection} on the tangent bundle 
		$TM$, produced by the COVID-19 dynamical system (\ref{1}), has the following form:%
		\[
		N=\left( N_{j}^{i}\right) _{i,j=\overline{1,6}}=\left( 
		\begin{array}{cccccc}
			0 & N_{2}^{1} & N_{3}^{1} & N_{4}^{1} & N_{5}^{1} & N_{6}^{1} \\ 
			- N_{2}^{1} & 0 & N_{3}^{2} & N_{4}^{2} & N_{5}^{2} & N_{6}^{2} \\  
			- N_{3}^{1} & -N_{3}^{2} & 0 & 0 & N_{5}^{3} & N_{6}^{3} \\ 
			- N_{4}^{1} & -N_{4}^{2} & 0 & 0 & 0 & N_{6}^{4} \\ 
			- N_{5}^{1} & -N_{5}^{2} & -N_{5}^{3} & 0 & 0 & N_{6}^{5} \\ 
			- N_{6}^{1} & -N_{6}^{2} & -N_{6}^{3} & -N_{6}^{4} & -N_{6}^{5} & 0 \\ 
		\end{array}%
		\right) , 
		\]%
		where%
		\[
		N_{2}^{1}=\dfrac{\beta_{s}I_{s}+\beta_{a}I_{a}+\beta_{h}I_{h}}{N^2} \cdot S +\dfrac{1}{2}\dfrac{\beta_{s}I_{s}+\beta_{a}I_{a}+\beta_{h}I_{h}}{N}, 
		\]
		\[
		N_{3}^{1}=-\dfrac{1}{2}\cdot\dfrac{-\beta_{s}N+\beta_{s}I_{s}+\beta_{a}I_{a}+\beta_{h}I_{h}}{N^2} \cdot S, 
		\]
		\[
		N_{4}^{1}=-\dfrac{1}{2}\cdot\dfrac{-\beta_{a}N+\beta_{s}I_{s}+\beta_{a}I_{a}+\beta_{h}I_{h}}{N^2} \cdot S,
		\]
		\[
		N_{5}^{1}=-\dfrac{1}{2}\cdot\dfrac{-\beta_{h}N+\beta_{s}I_{s}+\beta_{a}I_{a}+\beta_{h}I_{h}}{N^2} \cdot S,
		\]
		\[
		N_{6}^{1}=\dfrac{1}{2}\cdot\dfrac{\beta_{s}I_{s}+\beta_{a}I_{a}+\beta_{h}I_{h}}{N^2} \cdot S,
		\]
		\[
		N_{3}^{2}=-\dfrac{1}{2}\cdot \left [ \dfrac{\beta_{s}N-(\beta_{s}I_{s}+\beta_{a}I_{a}+\beta_{h}I_{h})}{N^2} \cdot S - (1-r)\sigma \right ],
		\]
		\[
		N_{4}^{2}=-\dfrac{1}{2}\cdot \left [ \dfrac{\beta_{a}N-(\beta_{s}I_{s}+\beta_{a}I_{a}+\beta_{h}I_{h})}{N^2} \cdot S - r\sigma \right ],
		\]
		\[
		N_{5}^{2}=-\dfrac{1}{2}\cdot \dfrac{\beta_{h}N-(\beta_{s}I_{s}+\beta_{a}I_{a}+\beta_{h}I_{h})}{N^2} \cdot S, 
		\]
		\[
		N_{6}^{2}=-\dfrac{1}{2}\cdot \dfrac{\beta_{s}I_{s}+\beta_{a}I_{a}+\beta_{h}I_{h}}{N^2} \cdot S , \;
		N_{5}^{3}=\dfrac{1}{2}\phi_{s}, \;
		N_{6}^{3}=\dfrac{1}{2}\gamma_{s}, \;
		\]
		\[
		N_{6}^{4}=\dfrac{1}{2}\gamma_{a}, \; 
		N_{6}^{5}=\dfrac{1}{2}\gamma_{h}.
		\]
	\end{proposition}
	
	Further, the Cartan canonical linear connection on the tangent bundle $TM$ produced by the
	least squares Lagrangian (\ref{4}) has all adapted components equal to
	zero. Moreover, the general formulas which give the Lagrangian d-torsions are expressed by (see 
	\cite{Mir-An})%
	\[
	R_{k}=\left( R_{jk}^{i}:={\frac{\delta N_{j}^{i}}{\delta x^{k}}}-\frac{%
		\delta N_{k}^{i}}{\delta x^{j}}\right) _{i,j=\overline{1,6}}, 
	\]%
	where%
	\[
	{\frac{\delta }{\delta x^{k}}={\frac{\partial }{\partial x^{k}}}-N_{k}^{r}{%
			\frac{\partial }{\partial y^{r}}.}} 
	\]%
	By direct computations, we obtain (see \cite{Bal-Nea})%
	\[
	R_{k}=\frac{\partial N}{\partial x^{k}},\quad \forall \;k=\overline{1,6}. 
	\]%
	Further, we infer
	
	\begin{proposition}
		The Lagrangian canonical Cartan linear connection, produced by the
		COVID-19 dynamical system (\ref{1}), is characterized by the $d-$\textbf{torsion skew-symmetric matrices}:%
		\[
		R_1=\dfrac{\partial N}{\partial S}, \; R_2=\dfrac{\partial N}{\partial E},\; R_3=\dfrac{\partial N}{\partial I_s}, \; R_4=\dfrac{\partial N}{\partial I_a},\; R_5=\dfrac{\partial N}{\partial I_h},\; R_6=\dfrac{\partial N}{\partial R}.
		\]
	\end{proposition}
	
	\begin{proposition}
		The\textbf{\ Lagrangian Yang-Mills e\-lec\-tro\-mag\-ne\-tic-like energy}%
		\textit{, produced by }the COVID-19 dynamical system (\ref{1}), is
		given by%
		\[
		\mathcal{EYM}(x)={\frac{1}{2}}\cdot Trace\left[ F\cdot F^{t}\right]= 
		\]
		\[
		=(N_{2}^{1})^2+(N_{3}^{1})^2+(N_{4}^{1})^2+(N_{5}^{1})^2+(N_{6}^{1})^2+(N_{4}^{2})^2+
		\]
		\[
		+(N_{5}^{2})^2+(N_{6}^{2})^2+(N_{5}^{3})^2+(N_{6}^{3})^2+(N_{6}^{4})^2+(N_{6}^{5})^2,
		\]
		where the electromagnetic-like matrix is $F=-N.$ For more details, see the works \cite{Mir-An} and \cite%
		{Bal-Nea}.
	\end{proposition}
	
	It is known that the \textit{matrix of deviation curvature} from the Kosambi-Cartan-
	Chern (KCC) geometrical theory is given by the formula
	\[
	P=\left( P_{j}^{i}\right) _{i,j=\overline{1,6}}=\frac{\partial N}{\partial
		x^{k}}y^{k}+\mathcal{E}, 
	\]%
	where, if%
	\[
	\mathcal{E}^{i}=2G^{i}-N_{j}^{i}y^{j}=-\frac{1}{2}\left( \frac{\partial X^{i}%
	}{\partial x^{j}}-\frac{\partial X^{j}}{\partial x^{i}}\right) y^{j}-\frac{%
		\partial X^{j}}{\partial x^{i}}X^{j} 
	\]%
	is the\textit{\ first invariant of the semispray} of the Lagrangian (\ref%
	{4}), then we put
	 $$\mathcal{E=}\left(\dfrac{\delta \mathcal{E}^{i}}{\delta x^{j}} \right) _{i,j=\overline{1,6}}.$$
	 Using direct calculations, we deduce that the entries of the matrix $\mathcal{E}$ are given by%
	\[
	\frac{\delta \mathcal{E}^{i}}{\delta x^{j}}=-\frac{1}{2}\left( \frac{%
		\partial ^{2}X^{i}}{\partial x^{j}\partial x^{k}}-\frac{\partial ^{2}X^{k}}{%
		\partial x^{i}\partial x^{j}}\right) y^{k}-\frac{\partial ^{2}X^{k}}{%
		\partial x^{i}\partial x^{j}}X^{k}-\frac{\partial X^{k}}{\partial x^{i}}%
	\frac{\partial X^{k}}{\partial x^{j}}- 
	\]%
	\[
	-\frac{1}{4}\left( \frac{\partial X^{k}}{\partial x^{j}}-\frac{\partial X^{j}%
	}{\partial x^{k}}\right) \left( \frac{\partial X^{i}}{\partial x^{k}}-\frac{%
		\partial X^{k}}{\partial x^{i}}\right) . 
	\]
	
	\begin{remark}
		Here above, all indices took values from $1$ to $6$. Moreover, the Einstein
		convention of summation was used.
	\end{remark}
	
	As a conclusion, it follows that the behavior of the neighboring solutions of the
	Euler-Lagrange equations ($\ref{5}$) is \textit{Jacobi stable} if and only if the real parts of the
	eigenvalues of the deviation tensor $P$ are strictly negative everywhere, and \textit{Jacobi unstable}, otherwise. The Jacobi stability or instability has the geometrical meaning
	that the trajectories of the Euler-Lagrange equations ($\ref{5}$) are bunching together or are dispersing (or are chaotic). For more details about KCC theory and its related
	Jacobi stability, see the works: B\"{o}hmer et al. $\cite{Bohmer}$, Bucătaru-Miron $\cite{Buc-Mir}$ and Neagu-Ovsiyuk $\cite{Nea-Ovs}$.
	
	\section{From COVID-19 dynamical system to Hamilton geometry}
	
	The \textit{least squares Hamiltonian }$H:T^{\ast
	}M\rightarrow \mathbb{R}$ associated with the Lagrangian (\ref{4}) is expressed by%
	\begin{equation} \label{6}
		\begin{array}{c}
			H(x,p)=\displaystyle{\frac{\delta ^{ij}}{4}}p_{i}p_{j}+X^{k}(x)p_{k}%
			\Leftrightarrow \medskip \\ 
			H(x,p)=\displaystyle{\frac{1}{4}}\left( p_{1}^{2}+p_{2}^{2}+p_{3}^{2}+p_{4}^{2}+p_{5}^{2}+p_{6}^{2}\right)+ \medskip \\ 
			+X^{1}(x)p_{1}+X^{2}(x)p_{2}+X^{3}(x)p_{3}+X^{4}(x)p_{4}+X^{5}(x)p_{5}+X^{6}(x)p_{6},%
		\end{array}
	\end{equation}%
	where $p_{r}=\partial L/\partial y^{r}$ and $H=p_{r}y^{r}-L$.
	
	Again, via the Hamilton geometry on cotangent bundles and the Hamiltonian least squares variational method for dynamical systems (see the monographs Miron et al. $\cite{Miron-Hr-Shim-Sab}$ and Neagu-Oană $\cite{Nea-Oana}$), it follows that we can develop a natural and distinct collection of nonzero Hamiltonian geometrical objects (such as nonlinear connection and d-torsions), which also characterize the
	COVID-19 dynamical system ($\ref{1}$).
	
	In this way, the Hamiltonian nonlinear connection on the cotangent bundle $T^{*}M$,
	has the components (see $\cite{Miron-Hr-Shim-Sab}$)
	\[
	N_{ij}=\frac{\partial ^{2}H}{\partial x^{j}\partial p_{i}}+\frac{\partial
		^{2}H}{\partial x^{i}\partial p_{j}}. 
	\]
	By direct computations, we find $\mathbf{N}=J+$ $J^{t}$. For more details,
	see the book $\cite{Nea-Oana}$. In conclusion, we get the following result:
	
	\begin{proposition}
		The \textbf{Hamiltonian canonical nonlinear connection} on the cotangent
		bundle $T^{\ast }M$, produced by the COVID-19 dynamical system (\ref%
		{1}), is described by the matrix $
		\mathbf{N}=\left( N_{ij}\right) _{i,j=\overline{1,6}}=$
		\[
		=\left( 
		\begin{array}{cccccc}
			N_{11} & N_{12} & N_{13} & N_{14} & N_{15} & N_{16}  \\ 
			N_{12} & N_{22} & N_{23} & N_{24} & N_{25} & N_{26}  \\  
			N_{13} & N_{23} & N_{33} & N_{34}=0 & N_{35}=\phi_{s} & N_{36}=\gamma_{s}  \\ 
			N_{14} & N_{24} & N_{43}=0 & N_{44}=-2\gamma_{a} & N_{45}=0 & N_{46}=\gamma_{a}  \\ 
			N_{15} & N_{25} & N_{53}=\phi_{s} & N_{54}=0 & N_{55} & N_{56}=\gamma_{h}  \\ 
			N_{16} & N_{26} & N_{63}=\gamma_{s} & N_{64}=\gamma_{a} & N_{56}=\gamma_{h} & N_{66}=0  \\ 
		\end{array}%
		\right) , 
		\]
	\end{proposition}	
	where 
	\[
	N_{11}=-2\dfrac{\beta_{s}I_{s}+\beta_{a}I_{a}+\beta_{h}I_{h}}{N^2} \cdot S -2\dfrac{\beta_{s}I_{s}+\beta_{a}I_{a}+\beta_{h}I_{h}}{N}, 
	\]
	\[
	N_{12}=\dfrac{\beta_{s}I_{s}+\beta_{a}I_{a}+\beta_{h}I_{h}}{N}, \;
	N_{13}=\dfrac{-\beta_{s}N+\beta_{s}I_{s}+\beta_{a}I_{a}+\beta_{h}I_{h}}{N^2} \cdot S,
	\]
	\[
	N_{14}=\dfrac{-\beta_{a}N+\beta_{s}I_{s}+\beta_{a}I_{a}+\beta_{h}I_{h}}{N^2} \cdot S,
	\]
	\[
	N_{15}=\dfrac{-\beta_{h}N+\beta_{s}I_{s}+\beta_{a}I_{a}+\beta_{h}I_{h}}{N^2} \cdot S,
	\]
	\[
	N_{16}=-\dfrac{\beta_{s}I_{s}+\beta_{a}I_{a}+\beta_{h}I_{h}}{N^2} \cdot S,
	\]
	\[
	N_{22}=2\left [ \dfrac{\beta_{s}I_{s}+\beta_{a}I_{a}+\beta_{h}I_{h}}{N^2} \cdot S -\sigma \right ],
	\]
	\[
	N_{23}=\dfrac{\beta_{s}N-(\beta_{s}I_{s}+\beta_{a}I_{a}+\beta_{h}I_{h})}{N^2} \cdot S +(1-r)\sigma,
	\]
	\[
	N_{24}=\dfrac{\beta_{a}N-(\beta_{s}I_{s}+\beta_{a}I_{a}+\beta_{h}I_{h})}{N^2} \cdot S +r\sigma,
	\]
	\[
	N_{25}=\dfrac{\beta_{h}N-(\beta_{s}I_{s}+\beta_{a}I_{a}+\beta_{h}I_{h})}{N^2} \cdot S,
	\]
	\[
	N_{26}=\dfrac{\beta_{s}I_{s}+\beta_{a}I_{a}+\beta_{h}I_{h}}{N^2} \cdot S, \;
	N_{33}=-2(\phi_{s}+\gamma_{s}+\delta_{s}), \;
	N_{55}=-2(\gamma_{h}+\delta_{h}).
	\]

	Moreover, the canonical Cartan linear connection on $T^{\ast }M$ produced by
	the least squares Hamiltonian (\ref{6}) has all adapted components equal
	to zero and the corresponding general formulas which give the Hamiltonian
	d-torsions are given by (see \cite{Miron-Hr-Shim-Sab})%
	\[
	\mathbf{R}_{k}=\left( R_{kij}:={\frac{\delta N_{ki}}{\delta x^{j}}}-{\frac{%
			\delta N_{kj}}{\delta x^{i}}}\right) _{i,j=\overline{1,6}}, 
	\]%
	where%
	\[
	{\frac{\delta }{\delta x^{j}}={\frac{\partial }{\partial x^{j}}}-N_{rj}{%
			\frac{\partial }{\partial p_{r}}.}} 
	\]%
	Consequently, we infer (for all details, see also \cite{Nea-Oana})
	
	\begin{proposition}
		The Hamiltonian canonical Cartan linear connection, produced by the
		COVID-19 dynamical system (\ref{1}), is characterized by the $d-$\textbf{torsion skew-symmetric matrices}:%
		\[
		\mathbf{R}_{k}=\frac{\partial }{\partial x^{k}}\left[ J-J^{t}\right]
		=-2R_{k},\quad \forall \;k=\overline{1,6}. 
		\]
	\end{proposition}
	
	\section{Conclusion}
	
	In this context, from our new geometric-physical approach, the hypersurfaces of constant level of the Lagrangian Yang-Mills electromagnetic-like energy produced by the COVID-19 dynamical system ($\ref{1}$) could have important connotations for the phenomenom taken in study. The hypersurface of constant level $\Sigma_\rho$ is obviously expressed by the equation
	
	\[
	(N_{2}^{1})^2+(N_{3}^{1})^2+(N_{4}^{1})^2+(N_{5}^{1})^2+(N_{6}^{1})^2+(N_{4}^{2})^2+
	\]
	\[
	+(N_{5}^{2})^2+(N_{6}^{2})^2+(N_{5}^{3})^2+(N_{6}^{3})^2+(N_{6}^{4})^2+(N_{6}^{5})^2=\rho\geq 0.
	\]
	
	The question is: \textit{There is a meaning, related to the spreading of the COVID-19 disease, for the shape of this hypersurface in the six dimensions $$(S,E,I_{s},I_{a},I_{h},R)=(x^1,x^2,x^3,x^4,x^5,x^6)?$$}
	
	For such a reason, we consider that the computer shown graphics of the surfaces determined by the 3-dimensional projections of the hypersurface $\Sigma_\rho$ on arbitrary three axes could offer important insights upon the spreading of the COVID-19 disease. This means that we talk about by the graphics of $C^3_6=20$ surfaces which could be relevant for the COVID-19 disease spreading.
	
	\noindent \textbf{Open problem.}
	What is the real meaning in epidemiology and spreading of the
	COVID-19 disease for our Lagrange-Hamilton geometrical objects constructed in
	this paper?

	\noindent \textsc{Ana-Maria Boldeanu and Mircea Neagu}\\[1mm]
	Transilvania University of Bra\c{s}ov\newline
	Department of Mathematics and Computer Science\newline
	Blvd. Iuliu Maniu, No. 50, 500091 Bra\c{s}ov, Romania.
	
	\noindent E-mails: ana.boldeanu@unitbv.ro, mircea.neagu@unitbv.ro
	
\end{document}